\documentclass[oneside,a4paper]{amsart} 
\usepackage{amssymb,verbatim,mathabx}

\usepackage[shortlabels]{enumitem}
    \setenumerate[0]{label={\rm (\roman*)}, leftmargin=*}

\usepackage[T1]{fontenc}
\usepackage[english]{babel}

\usepackage[textsize=footnotesize,textwidth=20ex,colorinlistoftodos]{todonotes}

\usepackage{hyperref}
\hypersetup{
    colorlinks=true,
    linkcolor=blue,
    filecolor=magenta, 
    citecolor=red,     
    urlcolor=cyan,
}

\usepackage[capitalize]{cleveref}
    \crefname{enumi}{}{}
    \Crefname{enumi}{Item}{Items}
    \crefname{equation}{}{}
    \Crefname{equation}{Equation}{Equations}

\newtheorem{proposition}{Proposition}[section]
\newtheorem{lemma}[proposition]{Lemma}

\newtheorem{theorem}[proposition]{Theorem}

\theoremstyle{definition}
\newtheorem{definition}[proposition]{Definition}

\newtheorem{example}[proposition]{Example}

\newtheorem{remark}[proposition]{Remark}

\begin{document}
\title{Line graphs of Multi-Graphs and the forbidden graph $E_6$}
\author{Hans Cuypers} 
\maketitle

\begin{abstract}
The \emph{line graph} $\Gamma$ of a multi-graph $\Delta$ 
is the graph whose vertices are the edges of $\Delta$,
where two such edges are adjacent if and only if they meet in a \emph{single} vertex of $\Delta$.

We provide several characterizations of such line graphs and in particular show that a graph is a line graph if and only if it does not contain one of $33$ graphs, all of which correspond to  bases
of anisotropic vectors of a $6$-dimensional  orthogonal geometry of $-$-type over a field with two elements, or, equivalently, to  sets of $6$ generating reflections in the Weyl group of type $E_6$.
\end{abstract}

\section{Introduction}
In this paper we consider (ordinary) graphs, i.e. graphs without loops and multiple edges, as well as 
multi-graphs, in which we allow multiple edges but no loops.
By a graph we will usually mean an ordinary graph.

Let $\Delta$ be a  multi-graph without loops but with possibly multiple edges. 
Then we define the \emph{line graph} of $\Delta$ to be the graph whose vertices are the edges of $\Delta$, where two edges are adjacent if and only if they meet in a single vertex.
We denote this line graph by $L(\Delta)$.
Notice that $L(\Delta)$ is an ordinary graph, a graph without loops or multiple edges.
(In \cite{cuypers_whitney}, the line graph of $\Delta$ is also called the $1$-line graph,
in contrast to the $\geq 1$-line graph of $\Delta$, which is the graph whose vertices
are the edges of $\Delta$, and two edges are adjacent if and only if there is at least one vertex on
both of them.)

\begin{figure}
{\tiny
\begin{tikzpicture}[scale=0.35]
\filldraw[color=black] 
(0,0) circle [radius=2pt]
(1,0) circle [radius=2pt]
(2,0) circle [radius=2pt]
(2,1) circle [radius=2pt]
(3,0) circle [radius=2pt]
(4,0) circle [radius=2pt];

\draw (0,0)--(4,0);
\draw (2,0)--(2,1);
\draw (2,-1.5) node {$E_6^{(1)}$};

\filldraw[color=black] 
(5,0) circle [radius=2pt]
(5,1) circle [radius=2pt]
(6,1) circle [radius=2pt]
(6,2) circle [radius=2pt]
(7,1) circle [radius=2pt]
(7,0) circle [radius=2pt];

\draw (7,0)--(5,0)--(5,1)--(7,1)--(7,0);
\draw (6,1)--(6,2);
\draw (6,-1.5) node {$E_6^{(2)}$};

\filldraw[color=black] 
(8,0) circle [radius=2pt]
(9,0) circle [radius=2pt]
(10,0) circle [radius=2pt]
(8,1) circle [radius=2pt]
(9,1) circle [radius=2pt]
(10,1) circle [radius=2pt];

\draw (8,0)--(10,0);
\draw (8,0)--(8,1)--(10,1);
\draw (9,0)--(9,1);
\draw (9,-1.5) node {$E_6^{(3)}$};

\filldraw[color=black] 
(11,0) circle [radius=2pt]
(12,0) circle [radius=2pt]
(13,0) circle [radius=2pt]
(14,0) circle [radius=2pt]
(12.5,1) circle [radius=2pt]
(13.5,1) circle [radius=2pt];

\draw (11,0)--(14,0);
\draw (13,0)--(12.5,1)--(13.5,1)--(13,0);
\draw (13,-1.5) node {$E_6^{(4)}$};

\filldraw[color=black] 
(16,0) circle [radius=2pt]
(17,0) circle [radius=2pt]
(18,0) circle [radius=2pt]
(16,1) circle [radius=2pt]
(17,1) circle [radius=2pt]
(18,1) circle [radius=2pt];

\draw (16,0)--(18,0)--(18,1)--(16,1)--(16,0);
\draw (17,0)--(17,1);
\draw (17,-1.5) node {$E_6^{(5)}$};

\filldraw[color=black] 
(19,0) circle [radius=2pt]
(20,0) circle [radius=2pt]
(19.5,2) circle [radius=2pt]
(20,1) circle [radius=2pt]
(21,0)  circle [radius=2pt]
(19,1) circle [radius=2pt];

\draw (19,0)--(20,0)--(20,1)--(19,1)--(19,0);
\draw (19,1)--(19.5,2)--(20,1);
\draw (20,0)--(21,0);
\draw (20,-1.5) node {$E_6^{(6)}$};

\filldraw[color=black] 
(22,0) circle [radius=2pt]
(23,0) circle [radius=2pt]
(22.5,2) circle [radius=2pt]
(23,1) circle [radius=2pt]
(24,1)  circle [radius=2pt]
(22,1) circle [radius=2pt];

\draw (22,0)--(23,0)--(23,1)--(22,1)--(22,0);
\draw (22,1)--(22.5,2)--(23,1);
\draw (23,1)--(24,1);
\draw (23,-1.5) node {$E_6^{(7)}$};

\filldraw[color=black] 
(25,1) circle [radius=2pt]
(26,1) circle [radius=2pt]
(26.5,2) circle [radius=2pt]
(26.5,0) circle [radius=2pt]
(27,1)  circle [radius=2pt]
(28,1) circle [radius=2pt];

\draw (25,1)--(26,1)--(26.5,2)--(26.5,0)--(27,1)--(28,1);
\draw (26,1)--(26.5,0);
\draw (26.5,2)--(27,1);
\draw (26.5,-1.5) node {$E_6^{(8)}$};

\filldraw[color=black] 
(29,0) circle [radius=2pt]
(29,1) circle [radius=2pt]
(30,0) circle [radius=2pt]
(31,0) circle [radius=2pt]
(31,1)  circle [radius=2pt]
(30,1) circle [radius=2pt];

\draw (29,0)--(31,0);
\draw (29,0)--(29,1)--(30,0)--(31,1)--(31,0);
\draw (30,0)--(30,1);
\draw (30,-1.5) node {$E_6^{(9)}$};

\filldraw[color=black] 
(32,0) circle [radius=2pt]
(33,0) circle [radius=2pt]
(34,0) circle [radius=2pt]
(35,0) circle [radius=2pt]
(33,1)  circle [radius=2pt]
(34,1) circle [radius=2pt];

\draw (32,0)--(35,0);
\draw (33,0)--(33,1)--(34,1)--(34,0)--(33,1);
\draw (34,-1.5) node {$E_6^{(10)}$};

\end{tikzpicture}

\medskip

\begin{tikzpicture}[scale=0.35]

\filldraw[color=black] 
(0,0) circle [radius=2pt]
(2,0) circle [radius=2pt]
(0,1) circle [radius=2pt]
(2,1) circle [radius=2pt]
(1,1)  circle [radius=2pt]
(1,2) circle [radius=2pt];

\draw (0,0)--(2,0)--(2,1)--(1,2)--(0,1)--(0,0);
\draw (0,0)--(1,1)--(2,0);
\draw (1,1)--(1,2);
\draw (1,-1.5) node {$E_6^{(11)}$};

\filldraw[color=black] 
(3,0) circle [radius=2pt]
(5,0) circle [radius=2pt]
(3,1) circle [radius=2pt]
(5,1) circle [radius=2pt]
(4,1)  circle [radius=2pt]
(4,2) circle [radius=2pt];

\draw (3,0)--(5,0)--(5,1)--(4,2)--(3,1)--(3,0);
\draw (4,1)--(4,2);
\draw (3,1)--(5,1);
\draw (4,-1.5) node {$E_6^{(12)}$};

\filldraw[color=black] 
(6,0) circle [radius=2pt]
(7,0) circle [radius=2pt]
(8,0) circle [radius=2pt]
(6,1) circle [radius=2pt]
(7,1)  circle [radius=2pt]
(8,1) circle [radius=2pt];

\draw (6,0)--(8,0)--(8,1)--(6,1)--(6,0)--(7,1)--(7,0);
\draw (7,-1.5) node {$E_6^{(13)}$};
\filldraw[color=black] 
(9,1) circle [radius=2pt]
(10,1) circle [radius=2pt]
(11,1) circle [radius=2pt]
(12,1) circle [radius=2pt]
(10,2)  circle [radius=2pt]
(10,0) circle [radius=2pt];

\draw (12,1)--(9,1)--(10,0)--(11,1);
\draw (10,0)--(10,1);
\draw (9,1)--(10,2)--(11,1);
\draw (10,-1.5) node {$E_6^{(14)}$};
\filldraw[color=black] 
(13,1) circle [radius=2pt]
(14,1) circle [radius=2pt]
(15,1) circle [radius=2pt]
(16,1) circle [radius=2pt]
(14,2)  circle [radius=2pt]
(14,0) circle [radius=2pt];

\draw (16,1)--(14,1);
\draw (14,0)--(14,2);
\draw (13,1)--(14,2)--(15,1);
\draw (13,1)--(14,0)--(15,1);
\draw (14,-1.5) node {$E_6^{(15)}$};
\filldraw[color=black] 
(17,1) circle [radius=2pt]
(18,1) circle [radius=2pt]
(19,1) circle [radius=2pt]
(17,0) circle [radius=2pt]
(19,0)  circle [radius=2pt]
(20,0) circle [radius=2pt];

\draw (17,1)--(19,1);
\draw (17,0)--(20,0);
\draw (17,0)--(17,1);
\draw (17,0)--(18,1);
\draw (19,0)--(18,1);
\draw (19,0)--(19,1);

\draw (18.5,-1.5) node {$E_6^{(16)}$};

\filldraw[color=black] 
(21,1) circle [radius=2pt]
(22,1) circle [radius=2pt]
(23,1) circle [radius=2pt]
(21,0) circle [radius=2pt]
(23,0)  circle [radius=2pt]
(22,2) circle [radius=2pt];

\draw (21,1)--(23,1);
\draw (21,0)--(23,0);
\draw (21,0)--(21,1);
\draw (21,0)--(22,1);
\draw (23,0)--(22,1);
\draw (22,1)--(22,2);
\draw (23,0)--(23,1);
\draw (22.5,-1.5) node {$E_6^{(17)}$};

\filldraw[color=black] 
(24,2) circle [radius=2pt]
(24,0) circle [radius=2pt]
(26,0) circle [radius=2pt]
(26,2) circle [radius=2pt]
(25,1)  circle [radius=2pt]
(27,1) circle [radius=2pt];

\draw (24,2)--(24,0)--(26,0)--(26,2)--(24,2);
\draw (24,2)--(26,0);
\draw (24,0)--(25,1);
\draw (26,2)--(27,1)--(26,0);

\draw (25,-1.5) node {$E_6^{(18)}$};

\filldraw[color=black] 
(28,2) circle [radius=2pt]
(28,0) circle [radius=2pt]
(30,0) circle [radius=2pt]
(30,2) circle [radius=2pt]
(29,1)  circle [radius=2pt]
(31,1) circle [radius=2pt];

\draw (28,2)--(28,0)--(30,0)--(30,2)--(28,2);
\draw (28,2)--(30,0);
\draw (30,2)--(29,1);
\draw (30,2)--(31,1)--(30,0);

\draw (29,-1.5) node {$E_6^{(19)}$};

\filldraw[color=black] 
(32,0.5) circle [radius=2pt]
(33,0) circle [radius=2pt]
(34,0) circle [radius=2pt]
(35,0.5) circle [radius=2pt]
(33,1)  circle [radius=2pt]
(34,1) circle [radius=2pt];

\draw (32,0.5)--(33,0)--(34,0)--(35,0.5)--(34,1)--(33,1)--(32,0.5);
\draw (33,0)--(34,1);
\draw (33,0)--(33,1);
\draw (34,0)--(34,1);
\draw (33.5,-1.5) node {$E_6^{(20)}$};

\end{tikzpicture}

\medskip

\begin{tikzpicture}[scale=0.35]
\filldraw[color=black] 
(-5,1) circle [radius=2pt]
(-4,0) circle [radius=2pt]
(-3,0) circle [radius=2pt]
(-2,1) circle [radius=2pt]
(-3,2)  circle [radius=2pt]
(-4,2) circle [radius=2pt];

\draw (-5,1)--(-4,0)--(-3,0)--(-2,1)--(-3,2)--(-4,2)--(-5,1);
\draw (-5,1)--(-3,0);
\draw (-5,1)--(-2,1);
\draw (-5,1)--(-3,2);
\draw (-5,1)--(-4,2);

\draw (-3.5,-1.5) node {$E_6^{(21)}$};

\filldraw[color=black] 
(0,0.5) circle [radius=2pt]
(1,0.5) circle [radius=2pt]
(2,0) circle [radius=2pt]
(2,1) circle [radius=2pt]
(3,0.5)  circle [radius=2pt]
(4,0.5) circle [radius=2pt];

\draw (0,0.5)--(1,0.5)--(2,0)--(2,1)--(1,0.5);
\draw (0,0.5)--(2,1)--(3,0.5)--(4,0.5);
\draw (0,0.5)--(2,0)--(3,0.5);

\draw (2,-1.5) node {$E_6^{(22)}$};

\filldraw[color=black] 
(5,0.5) circle [radius=2pt]
(6,0.5) circle [radius=2pt]
(7,0) circle [radius=2pt]
(7,1) circle [radius=2pt]
(8,0.5)  circle [radius=2pt]
(7,2) circle [radius=2pt];

\draw (5,0.5)--(6,0.5)--(7,0)--(7,1)--(6,0.5);
\draw (5,0.5)--(7,1)--(8,0.5);
\draw (7,1)--(7,2);
\draw (5,0.5)--(7,0)--(8,0.5);

\draw (7,-1.5) node {$E_6^{(23)}$};
\filldraw[color=black] 
(9,1) circle [radius=2pt]
(10,0) circle [radius=2pt]
(11,0) circle [radius=2pt]
(12,1) circle [radius=2pt]
(11,2)  circle [radius=2pt]
(10,2) circle [radius=2pt];

\draw (9,1)--(10,0)--(11,0)--(12,1)--(11,2)--(10,2)--(9,1);
\draw (10,0)--(10,2);
\draw (11,0)--(11,2);
\draw (9,1)--(11,0);
\draw (12,1)--(10,0);

\draw (10.5,-1.5) node {$E_6^{(24)}$};
\filldraw[color=black] 
(13,1) circle [radius=2pt]
(14,0) circle [radius=2pt]
(15,0) circle [radius=2pt]
(16,1) circle [radius=2pt]
(14.5,2)  circle [radius=2pt]
(14.5,1) circle [radius=2pt];

\draw (13,1)--(14,0)--(15,0)--(16,1)--(14.5,2)--(13,1)--(14.5,1)--(14,0);
\draw (15,0)--(14.5,1)--(16,1);
\draw (14.5,2)--(14.5,1);

\draw (14.5,-1.5) node {$E_6^{(25)}$};

\filldraw[color=black] 
(17,0) circle [radius=2pt]
(17,2) circle [radius=2pt]
(18,1) circle [radius=2pt]
(19,1) circle [radius=2pt]
(20,1)  circle [radius=2pt]
(21,1) circle [radius=2pt];

\draw (17,0)--(17,2)--(21,1)--(18,1)--(17,0)--(18,1)--(17,2);
\draw (21,1)--(17,0)--(19,1)--(17,2);

\draw (18,-1.5) node {$E_6^{(26)}$};

\filldraw[color=black] 
(22,1) circle [radius=2pt]
(23,0) circle [radius=2pt]
(24,0) circle [radius=2pt]
(25,1) circle [radius=2pt]
(23.5,2)  circle [radius=2pt]
(23.5,1) circle [radius=2pt];

\draw (22,1)--(23,0)--(24,0)--(25,1)--(23.5,2)--(22,1);

\draw (23.5,2)--(24,0)--(23.5,1);
\draw (23.5,2)--(23,0)--(23.5,1);

\draw (23.5,2)--(23.5,1);

\draw (23.5,-1.5) node {$E_6^{(27)}$};

\end{tikzpicture}

\medskip

\begin{tikzpicture}[scale=0.35]

\filldraw[color=black] 
(-5,0) circle [radius=2pt]
(-5,2) circle [radius=2pt]
(-4,1) circle [radius=2pt]
(-3,1) circle [radius=2pt]
(-2,1)  circle [radius=2pt]
(-1,1) circle [radius=2pt];

\draw (-5,0)--(-5,2)--(-4,1)--(-5,0);
\draw (-4,1)--(-1,1);
\draw (-5,0)--(-4,1)--(-5,2);
\draw (-5,0)--(-3,1)--(-5,2);
\draw (-5,0)--(-2,1)--(-5,2);

\draw (-4,-1.5) node {$E_6^{(28)}$};

\filldraw[color=black] 
(0,1) circle [radius=2pt]
(1,0) circle [radius=2pt]
(2,0) circle [radius=2pt]
(3,1) circle [radius=2pt]
(2,2)  circle [radius=2pt]
(1,2) circle [radius=2pt];

\draw (0,1)--(1,0)--(2,0)--(3,1)--(2,2)--(1,2)--(0,1);
\draw (1,0)--(1,2);
\draw (2,0)--(2,2);
\draw (1,2)--(2,0);
\draw (2,2)--(1,0);
\draw (0,1)--(3,1);
\draw (1.5,-1.5) node {$E_6^{(29)}$};

\filldraw[color=black] 
(4,1) circle [radius=2pt]
(5,1) circle [radius=2pt]
(6,2) circle [radius=2pt]
(6,0) circle [radius=2pt]
(7,1)  circle [radius=2pt]
(8,1) circle [radius=2pt];

\draw (4,1)--(6,2)--(8,1)--(6,0);
\draw (4,1)--(5,1)--(6,2)--(7,1)--(8,1);
\draw (5,1)--(6,0)--(7,1);
\draw (4,1)--(6,0)--(6,2);

\draw (6,-1.5) node {$E_6^{(30)}$};

\filldraw[color=black] 
(9,1) circle [radius=2pt]
(10,2) circle [radius=2pt]
(10,0) circle [radius=2pt]
(11,0.7) circle [radius=2pt]
(11,1.3)  circle [radius=2pt]
(12.5,1) circle [radius=2pt];

\draw (9,1)--(10,2)--(12.5,1)--(10,0)--(9,1);
\draw (10,0)--(10,2)--(11,1.3)--(11,0.7)--(10,0)--(11,1.3);
\draw (11,1.3)--(12.5,1)--(11,0.7);

\draw (10,-1.5) node {$E_6^{(31)}$};

\filldraw[color=black] 
(14,1) circle [radius=2pt]
(15,1) circle [radius=2pt]
(16,2) circle [radius=2pt]
(16,0) circle [radius=2pt]
(17,1)  circle [radius=2pt]
(18,1) circle [radius=2pt];

\draw (14,1)--(16,2)--(18,1)--(16,0);
\draw (14,1)--(15,1)--(16,2)--(17,1)--(18,1);
\draw (15,1)--(16,0)--(17,1);
\draw (14,1)--(16,0)--(16,2);
\draw (15,1)--(17,1);
\draw (16,-1.5) node {$E_6^{(32)}$};

\filldraw[color=black] 
(19,1) circle [radius=2pt]
(20,0) circle [radius=2pt]
(21,0) circle [radius=2pt]
(22,1) circle [radius=2pt]
(21,2)  circle [radius=2pt]
(20,2) circle [radius=2pt];

\draw (19,1)--(20,0)--(21,0)--(22,1)--(21,2)--(20,2)--(19,1);
\draw (19,1)--(22,1);
\draw (20,0)--(21,2);
\draw (21,0)--(20,2);
\draw (20,0)--(20,2);

\draw (19,1)--(21,0);
\draw (22,1)--(20,2);
\draw (20,0)--(22,1);

\draw (20.5,-1.5) node {$E_6^{(33)}$};
\end{tikzpicture}
}

\caption{The graphs equivalent to $E_6$.}
\label{E6graphs}
\end{figure}
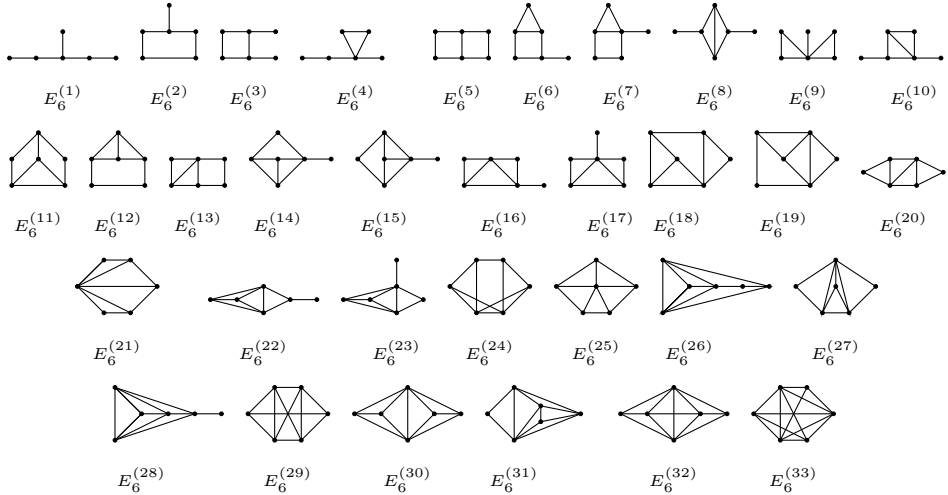

We prove the following:

\begin{theorem}\label{mainthm}
Let $\Gamma$ be a  connected ordinary graph.
Then $\Gamma$ is a line graph of a multi-graph
if and only if it does not contain an induced subgraph of the set $\mathcal{E}_6$
of $33$ graphs given  in \cref{E6graphs}.
\end{theorem}

The graphs in $\mathcal{E}_6$   can be described as follows.
Consider the orthogonal space $(V,Q)$ of dimension $6$ over $\mathbb{F}_2$ and $Q$ a quadratic form of $-$-type (i.e. with Arf invariant $+1$). 
For any basis of  this space consisting of anisotropic vectors (vectors $v$ with $Q(v)=1$) we consider   
the graph whose vertices are the vectors in the basis and two vertices $v\neq w$ being adjacent if and only if $Q(v+w)=1$.
The set  $\mathcal{E}_6$ consists of all such connected graphs and equals the set of $33$ graphs given in \cref{E6graphs}. The most prominent graph in this set is of course the graph $E_6$.
Alternatively, one can describe the graphs in $\mathcal{E}_6$ as the graphs on $6$ generating reflections in the Weyl group of type $E_6$, where 
two reflections are adjacent if and only if they do not commute, or as the graph on $6$ roots spanning
the root lattice $E_6$, where two such roots are adjacent if and only if they are not perpendicular. 

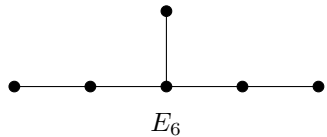
\begin{center}
\begin{tikzpicture}
\filldraw[color=black] 
(0,0) circle [radius=2pt]
(1,0) circle [radius=2pt]
(2,0) circle [radius=2pt]
(2,1) circle [radius=2pt]
(3,0) circle [radius=2pt]
(4,0) circle [radius=2pt];

\draw (0,0)--(4,0);
\draw (2,0)--(2,1);
\draw (2,-0.5) node {$E_6$};
\end{tikzpicture}

\end{center}
We call a graph \emph{$\mathcal{E}_6$-free}, if it does not contain an induced subgraph from $\mathcal{E}_6$.

\cref{mainthm} is closely related to the 
following characterizations of  line graphs of ordinary graphs
by  nine forbidden induced subgraphs  due to Beineke \cite{beineke} and of generalized line graphs
by $31$ forbidden induced  subgraphs by Cvetkovi\'{c}, Doob and Simi\'{c} \cite{GLG} for finite graphs and Vijayakumar \cite{vija2} for infinite graphs.

\begin{theorem}[Beineke, \cite{beineke}]\label{LGthm}
Let $\Gamma$ be a connected graph.
Then $\Gamma$ is a line graph of an ordinary graph
if and only if it does not contain one of $9$ graphs  
$H_1,H_2, H_3$ as in \cref{forbiddenlinegraphs},
or of the graphs $E_6^{(8)},E_6^{(12)},E_6^{(20)},E_6^{(22)}, E_6^{(25)}$ or $E_6^{(30)}$ from \cref{E6graphs}.
\end{theorem}

\begin{theorem}[Cvetkovi\'{c} et al., \cite{GLG}, Vijayakumar, \cite{vija2}]\label{GLGthm}
Let $\Gamma$ be a connected graph.
Then $\Gamma$ is a generalized line   graph if and only if $\Gamma$ does not contain one of $11$ graphs $G_1,\dots, G_{11}$ from \cref{Ggraphs}
or one of the $20$ graphs $E_6^{(i)}$, where $i\in \{1,2,4,8,9,10,12,16,17,20,21,23,24,25,27,28,30,31,32,33\}$,   from \cref{E6graphs}. 
\end{theorem}

Actually, both these results follow easily from \cref{mainthm}.

\cref{mainthm} not only generalizes the above two results, but also provides a common explanation why the minimal forbidden graphs in each of the three theorems above have at most $6$ vertices by relating them to the Weyl group and  root lattice of type $E_6$.

There where the main tool for studying generalized line graphs, or more generally graphs with least eigenvalue $\geq -2$, is embedding them into a root lattice, see \cite{css},  
our main tool in proving the above theorems is to consider the 
graphs as embedded into an orthogonal space over $\mathbb{F}_2$,
the field of $2$ elements. 
We notice that the  $E_6$-root lattice modulo $2$ corresponds to the $6$-dimensional 
orthogonal space of $-$-type, whose bases provide us with the class $\mathcal{E}_6$.

As our approach relates to the theory of Coxeter and Weyl groups, root lattices and related geometries and groups,
see \cite{brown,brown2,css,hall,seven}, we adopt the names of Dynkin diagrams $A_n$, $D_n$ and $E_n$ for the corresponding graphs.

The remainder of this paper is organised as follows.

For any ordinary graph $\Gamma$ we can view its vertices
as certain anisotropic vectors in an orthogonal space $(V,Q)$,
where two vertices are adjacent if and only if $Q(v+w)=1$.
We can of course assume that the  vertices of $\Gamma$ do linearly span $V$, but that does not imply that they generate the 
partial linear space of anisotropic vectors and elliptic lines of the orthogonal space.
The case that the vertices of $\Gamma$ do not generate this partial linear space can be shown to imply that $\Gamma$ is the line graph of a  multi-graph.
The embedding of a graph $\Gamma$ into an orthogonal space over $\mathbb{F}_2$ is the topic of \cref{sect:embedding}.

To check whether the points of a graph $\Gamma$ embedded in an orthogonal $\mathbb{F}_2$-space $(V,Q)$ generate the full 
partial linear space of anisotropic vectors and elliptic lines, we can change the graph $\Gamma$ into a new graph $\Delta$ generating
the same subspace.
This provides us with an equivalence relation between graphs, which is studied in \cref{sect:equivalence} and \cref{sect:E6} and leads to a proof of \cref{mainthm}. 

In these sections we make use of techniques and results from 
Brown and Humphreys  \cite{brown, brown2}
and Seven \cite{seven}.
We show that each equivalence class of graphs  contains  a tree.
In particular, we reprove a remarkable result by Seven
(Seven's \cref{seven}) that any graph $\Gamma$ equivalent to a tree with an $E_6$ subgraph contains itself a subgraph from the equivalence class $\mathcal{E}_6$ of the graph $E_6$.
This lemma is the key to prove \cref{mainthm}.

For completeness we provide (partly new and shorter) proofs of the results of \cite{brown, brown2,seven} needed in this paper, based on 
a result presented by Arjeh Cohen at his PhD defence and part of the Appendix of \cite{cuypers}.

In the  two sections  \ref{sect:linegraph} and \ref{sect:genlinegraph} we prove the two theorems
\cref{LGthm} and \cref{GLGthm}.
We end the paper with some final remarks in \cref{sect:remarks}.

\medskip

\noindent
{\bf Acknowledgement}.
The author thanks Jonathan I. Hall for bringing Seven's lemma
to his attention, and Arjeh Cohen for sharing the results of the Appendix of \cite{cuypers}.

\section{The orthogonal embedding of a graph}
\label{sect:embedding}

Let $\Gamma=(\mathcal{V},\mathcal{E})$ be a  graph 
with vertex set $\mathcal{V}$ and edge set $\mathcal{E}$. 
Then consider ${V}_\Gamma$, the vector space of finite
subsets of $\mathcal{V}$, where for two finite subsets $u,w$ of $\mathcal{V}$ the sum
$u+w$ is defined to be the symmetric difference of $u$ and $w$.
We identify a vertex $v$ with the subset $\{v\}\in V_\Gamma$. 

Put a total ordering $<$ on the vertex set of $\Gamma$.
Let $u$ and $w$ be two finite subsets of $\mathcal{V}$ and let $g_\Gamma(u,w)$ denote the number of ordered pairs $(x,y)\in u\times w$, where $x<y$ and $\{x,y\}$ is an edge, or $x=y$, modulo $2$.
Then  $g_\Gamma(u,v+w)=g_\Gamma(u,v)+g_\Gamma(u,w)$ for any finite subsets $u,v,w$ of $\mathcal{V}$,
as the ordered edges $(x,z)$ with $x<z$ and $z\in v\cap w$ are counted twice at the right hand site of the equation, just as vertices in the intersection of $u$ and $v\cap w$.

Similarly we find $g_\Gamma(v+ w,u)=g_\Gamma(v,u)+g_\Gamma(w,u)$. So,  $g_\Gamma:{V}_\Gamma\times {V}_\Gamma\rightarrow \mathbb{F}_2$ is bilinear.
The map $Q_\Gamma:V\rightarrow \mathbb{F}_2$ given by $Q_\Gamma(v)=g_\Gamma(v,v)$ for all $v\in V_\Gamma$ is a quadratic form with associated symmetric (and also alternating) form
$f_\Gamma$ given by $f_\Gamma(u,w)=g_\Gamma(u,w)+g_\Gamma(w,u)$.

Notice that if we take a different total ordening $<'$, and let $g_\Gamma'$ be the corresponding bilinear map,
then $g_\Gamma(u,w)+g_\Gamma(w,u)=g_\Gamma'(u,w)+g_\Gamma'(w,u)$ and $g_\Gamma(v,v)=g_\Gamma'(v,v)$. So, $Q_\Gamma$ and $f_\Gamma$ are independent of the chosen ordening.

An \emph{orthogonal embedding} of a graph $\Gamma=(\mathcal{V},\mathcal{E})$ is a map from $\mathcal{V}$ into the set of anisotropic vectors of an orthogonal $\mathbb{F}_2$-space  $(V,Q)$,
such that the images of the vertices span $V$ and  two vertices $v$ and $w$ of $\Gamma$ are adjacent if and only if the sum of their images is also anisotropic.
The embedding of $\Gamma$ into $(V_\Gamma,Q_\Gamma)$, mapping a vertex $v$ to $\{v\}$, is an orthogonal embedding.
It is called the \emph{universal orthogonal embedding} of $\Gamma$.

Now suppose $\mathcal{V}$ is a subset of anisotropic points of an orthogonal $\mathbb{F}_2$-space $(V,Q)$, then the graph with vertex set $\mathcal{V}$ in which two vertices $v,w$
are adjacent if and only if $Q(v+w)=1$ is called the graph \emph{induced} on $\mathcal{V}$ by $Q$. 
So, a graph $\Gamma$ is (isomorphic to) the induced graph of its universal orthogonal embedding into $(V_\Gamma,Q_\Gamma)$.

Let $\Gamma$ be a graph and  $R$  a subspace of the isotropic radical $\{v\in V_\Gamma\mid Q_\Gamma(v)=0=f_\Gamma(v,w)$ for all $w\in V_\Gamma\}$ of $(V_\Gamma,Q_\Gamma)$, then we can take the quotient modulo $R$ and
find an orthogonal embedding of $\Gamma$ in the quotient space. If we take for $R$ the full isotropic radical, then we call the  embedding of $\Gamma$ into $V_\Gamma/R$ the 
\emph{minimal (orthogonal) embedding} of $\Gamma$.
We denote the graph induced on the images of the elements of $\mathcal{V}$ in the minimal embedding of $\Gamma$
by $\overline{\Gamma}=(\overline{\mathcal V},\overline{\mathcal E})$.

Obviously we have:

\begin{proposition}\label{orth}
Let $\Gamma$ be a connected graph.
\begin{enumerate}
\item Every orthogonal embedding of $\Gamma$ can be obtained as a quotient
of the universal embedding by some subspace contained in the isotropic radical of  $(V_\Gamma,Q_\Gamma)$.
\item In the minimal orthogonal embedding of $\Gamma$  two vertices are mapped to the same vector if and only if they are non-adjacent and have the same set of neighbors.
\item $\Gamma$ is isomorphic to $\overline{\Gamma}$ if and only if any two vertices $v,w\in \mathcal{V}$ do not have the same set of neighbors.

\end{enumerate}
\end{proposition}

Let $\Gamma$ be a connected graph, then two adjacent vertices $v,w$ are mapped to 
two vectors in $V_\Gamma$, also denoted by $v$ and $w$, such that $$ Q_\Gamma(v)=Q_\Gamma( w)= Q_\Gamma(v+ w)=1.$$
This implies that the images of the vertices of $\Gamma$ generate
a connected subspace of the \emph{orthogonal cotriangular space} of $(V_\Gamma,Q_\Gamma)$, i.e.  the partial linear space of anisotropic points outside the radical
and elliptic lines of $(V_\Gamma,Q_\Gamma)$. Here an \emph{elliptic line} is considered to be a set of three anisotropic vectors in a $2$-dimensional subspace on which $Q_\Gamma$ takes the value $1$ for each non-zero vector.
We denote this subspace by $\Pi_\Gamma$. 
This space  maps onto a subspace the cotriangular space  in the quotient $(\overline V_\Gamma, \overline Q_\Gamma)$,
which we then denote by $\overline\Pi_\Gamma$.

\begin{lemma} 
$\overline \Pi_\Gamma=\Pi_{\overline\Gamma}$.
\end{lemma}

\begin{proof}
The vertices of $\overline\Gamma$ are contained in $\overline \Pi_\Gamma$,
so $\Pi_{\overline\Gamma}$ is a subspace of $\overline \Pi_\Gamma$.

Now suppose $\overline v$ is a point of  $\overline \Pi_\Gamma$.
Then there are vertices $v_1,\dots, v_k$ of $\Gamma$ with  $v_1+\dots+v_i$ adjacent to $v_{i+1}$ such that $v=v_1+\cdots+ v_k$ is a point of $\Pi_\Gamma$
mapping to $\overline v$.
But then also $\overline v=\overline v_1+\cdots +\overline v_k$ with  $\overline{v_1+\dots+v_i}=
\overline{v_1}+\dots+\overline{v_i}$ adjacent to $\overline{v_{i+1}}$ in $\overline\Gamma$, and we find $\overline v$ to be a point
of $\Pi_{\overline\Gamma}$.
\end{proof}

\begin{example}\label{Anexample}
Let $\Omega$ be a set and 
$\Gamma$ the line graph of the complete graph $\Delta$ on $\Omega$ with vertices $v_{ij}:=\{i,j\}$, where $i\neq j\in \Omega$.

Then the isotropic radical of $(V_\Gamma,Q_\Gamma)$
contains $v_{ij}+v_{jk}+v_{kl}$, where $\{i,j,k\}$ is a subset of size $3$ of $\Omega$.

The partial linear  space $\overline\Pi_\Gamma$ is 
generated by the subset of images modulo the isotropic radical of vertices $v_{ij}$ of any line graph of a connected
graph $\Delta_0$ on the vertex set $\Omega$.

Indeed, if $\overline v_{ij}$ is a point of $\overline\Pi_\Gamma$, with $v_{ij}$  not an edge of $\Delta_0$,
then there is a finite path in $\Delta_0$ from $i$ to $j$.
Let $v_{ik_1},v_{k_1k_2},\dots, v_{k_n j}$ be the edges involved in this path.
Then $\overline v_{ij}=\overline v_{ik_1}+\overline v_{k_1k_2}+\cdots+ \overline v_{k_n j}$
is clearly in  $\overline\Pi_\Gamma$.

The space  $\overline\Pi_\Gamma$  is isomorphic to the \emph{cotriangular space of the set} $\Omega$, i.e. the partial linear space with as points the pairs from $\Omega$ and as lines the triples of pairs in subsets of size $3$, unless $|\Omega|=4$, in which case  $\overline\Pi_\Gamma$ is isomorphic to
the cotriangular space of a set of size $3$.

If $\Omega=\{1,\dots, 6\}$, then  $(\overline V_\Gamma,\overline Q_\Gamma)$ is a
$5$-dimensional  orthogonal
space with a trivial isotropic radical and $\Pi_\Gamma$ coincides with the orthogonal cotriangular space on $15$ points.
If $\Omega=\{1,\dots, 8\}$, then
the vector $v_{12}+v_{34}+v_{56}+v_{78}$ is in the isotropic radical of $Q_\Gamma$, and 
$(\overline V_\Gamma,\overline Q_\Gamma)$ is a  $6$-dimensional nondegenerate orthogonal
space of $+$-type. Again $\Pi_\Gamma$ is the full orthogonal cotriangular space  on all $28$ anisotropic vectors and elliptic lines of $ \overline V_\Gamma$.
\end{example}

\begin{example}\label{E6example}
If $\Gamma$ is the graph $E_6$, then $(V_{\Gamma},Q_{\Gamma})$ is a nondegenerate $6$-dimensional space of $-$-type.
The partial linear space $\Pi_{\Gamma}$ contains all $36$ anisotropic vectors.
We find
 $\Pi_{\Gamma}$ not to be isomorphic to  the cotriangular space of some set.
A way of seeing this, which we will use later, is that $\Pi_{\Gamma}$ contains
a subspace $S$ isomorphic to $\Pi_{D_4}$ and consisting of $12$ points.
The subspace $S$ generates a $4$-dimensional subspace of
$V_{\Gamma}$ with a $2$-dimensional radical.
Indeed, the subspace generated by the $4$ vertices of the $D_4$ subgraph of $\Gamma$ provides us with such a subspace.
Such subspaces do not exist in any cotriangular space of a set.

We can phrase the above also in terms of subgroups of the orthogonal group of $(V_{\Gamma},Q_{\Gamma})$.

Each anisotropic vector $v$ of $V_{\Gamma}$
is the center of a transvection
defined by $$w\in V_{\Gamma}\mapsto w+f_{\Gamma}(v,w)v.$$
The $6$ transvections associated to the vertices 
of $\Gamma$ generate the orthogonal group
$O(V_{\Gamma},Q_{\Gamma})$ which is isomorphic to
the Weyl group of type $E_6$.

Inside this group the $4$ transvections associated to
the subgraph of type $D_4$ of $\Gamma_0$ generate
a subgroup isomorphic to $W(D_4)$, the Weyl group type $D_4$.

If we consider the minimal orthogonal embedding 
of the line graph $\Gamma'$ of a complete graph on $n\geq 5$
vertices, then 
the transvections corresponding to the vertices of 
$\Gamma'$
generate a group isomorphic to the symmetric group on $n$ letters, i.e., a Weyl group of type $A_{n-1}$. 
Inside this group, there are no $4$ transvections generating a group isomorphic to $W(D_4)$.
\end{example}

\begin{theorem}\label{plsthm}
Let $\Gamma$ be a connected ordinary graph.
Then the following statements are equivalent:
\begin{enumerate}
\item $\Gamma$ is the line graph of some multi-graph.
\item $\overline\Gamma$ is the line graph of some  ordinary graph.
\item $\overline\Pi_{\Gamma}$ is isomorphic to a cotriangular space of some set $\Omega$.
\item $\overline\Pi_{\Gamma}$ contains no subspace isomorphic to $\Pi_{D_4}$.
\end{enumerate}
\end{theorem}

\begin{proof}
We first show that (i) implies (ii). Assume $\Gamma$ is the line graph of a connected multi-graph $\Delta$. Then two vertices $v$ and $w$ corresponding to edges on the same vertices in $\Delta$ have the same set of neighbors in $\Gamma$.
This implies that they are identified in $\overline\Gamma$.
So, if $\underline \Delta$ is the graph obtained from $\Delta$ by replacing all multiple edges
with single edges, we find that $\overline{L(\underline\Delta)}=\overline\Gamma$.

If $\overline{L(\underline\Delta)}= L(\underline\Delta)$, we are done.
So, assume that $L(\underline\Delta)$ contains two vertices $v$ and $w$ that are nonadjacent
but have the same set of neighbors.
That implies that $\underline\Delta$ is a connected  graph on $4$ vertices, and 
$\overline\Gamma$ a graph on at most $3$ vertices, which is clearly a line graph of an ordinary graph.
So (i) implies (ii).

 Now assume $\overline\Gamma$ is a line graph of an ordinary graph $\underline \Delta$.
Then $\Gamma$ can be seen to be the line graph of the graph $\Delta$ with the same vertex set as $\underline\Delta$, but in which the edges $e_{\overline v}$ of $\underline \Delta$ corresponding to  vertices $\overline v$ of $\overline \Gamma$ are replaced by multiple edges $e_v$ on the same vertices in $\Delta$,
where $v$ runs over the vertices in $\Gamma$ that are mapped to $\overline v$.
So (ii) implies (i).

Let $\Gamma$ be the line graph of a multi-graph. Then, by the above we find $\overline \Gamma$ to be a line graph of an ordinary graph. But then \cref{Anexample} shows that  $\overline\Pi_{\Gamma}=\Pi_{\overline\Gamma}$ is  a cotriangular space of some set $\Omega$. So (ii) implies (iii).

Now suppose $\overline\Pi_{\Gamma}=\Pi_{\overline\Gamma}$ equals a cotriangular space  $\Pi$
of a set $\Omega$.

Then all the vertices of $\overline\Gamma$ can be identified with pairs from
$\Omega$. In particular,  $\overline\Gamma$ is a line graph. This shows that (iii) implies (ii).

Moreover, inside $\Pi_{\overline\Gamma}$ we do not find subspaces isomorphic to $\Pi_{D_4}$.
So (iii) also implies (iv).

Now assume  $\overline\Pi_{\Gamma}$ contains no subspace isomorphic to $\Pi_{D_4}$.
Then let $v,w$ be two collinear points and $u=v+w$ the third point on the line on $v$ and $w$.
Any point $x$ is collinear with $0$ or two points of $\{u,v,w\}$.
Moreover, if $x$ and $y$ are two points collinear with both $v,w$, but not with $u$,
then $x$ and $y$ are collinear, for otherwise we find a subgraph $D_4$ on 
$v,u,x,y$. So, $\omega_{v,w}$, defined as the set of points $x$ of $\overline\Pi_{\Gamma}$ that
are collinear to $v,w$ but not $u$, is a clique. 

If $x$ is a point in $\omega_{v,w}$ different from $v,w$, then any point
$y\in \omega_{v,w}$  different from $x,v$ and $w$,  
is collinear with $x,v$ and $w$, and then different from $v+x$.
So, $\omega_{v,w}\subseteq \omega_{x,v}$ and by symmetry of the argument  $\omega_{v,w}\subseteq \omega_{x,v}$.
But the we find $\omega_{v,w}=\omega_{x,y}$ for any two points in $\omega_{v,w}$.

Let $\Omega$ be the set of all cliques $\omega_{x,y}$, where $x,y$ are collinear points of $\overline\Pi_{\Gamma}$.

The point $v$ is contained in both $\omega_{v,w}$ and in $\omega_{v,u}$.
Let $\omega$ be an arbitrary element from $\Omega$  containing $v$.
So, $\omega=\omega_{v,x}$ for some point $x$ collinear to $x$.
If $x$ is different from $w,u$, then it is collinear with one of  $w$ or $u$
and we find $w$ or $u$ is in $\omega$ from which we deduce that
$\omega$ equals $\omega_{v,w}$ or $\omega_{v,u}$.
So $v$ is in exactly two elements of $\Omega$ and two elements from $\Omega$ intersect on at most one point.
But then we can identify the points of $\overline\Gamma$ with pairs of elements from $\Omega$
and find $\overline\Gamma$ to be a line graph. So (iv) implies (ii).
\end{proof}

\begin{remark}
The above proof is based upon the ideas of Krausz \cite{krausz} characterizing line graphs of ordinary graphs,  and of  Hall \cite[Section 6]{general}, characterizing the (finitary) symmetric groups as $3$-transposition groups.
\end{remark}

\section{An equivalence between Graphs}
\label{sect:equivalence}

Let $\Gamma=(\mathcal{V},\mathcal{E})$ be a  graph with vertex set $\mathcal{V}$ and edge set
$\mathcal{E}$. 
Without loss of generality we can assume $\Gamma$ to be the induced subgraph on the set $\mathcal{V}$ of anisotropic
vectors in an orthogonal $\mathbb{F}_2$-space $(V,Q)$. Denote by $f$ the symmetric bilinear form associated to $Q$.

For vertices $v,w\in \mathcal{V}$ we define  
$\Gamma^{(v,w)}$ to be the graph induced on $\mathcal{V}\setminus\{w\}\cup \{w+f(v,w)v\}$.
So, we replace the vertex $w$ by $\tau_v(w):=w+f(v,w)v$, all other vertices remain the same.

Such transformation is called an elementary transformation.
We say that $\Gamma$ is \emph{equivalent} to
$\Delta$ if it can be obtained from $\Delta$ by a series of elementary
transformations.

Clearly, equivalence is an equivalence relation.
This relation has been studied by Brown and Humphries \cite{brown,brown2} as well as Seven \cite{seven}.

In this section we analyse the equivalence classes of graphs.
We use the approach of the  Appendix in \cite{cuypers}.

\begin{lemma}\label{connected}
Let $v,w$ be two vertices of an ordinary graph $\Gamma$.
\begin{enumerate}
\item
If $\Gamma$ is connected, then so is $\Gamma^{(v,w)}$.
\item $\Pi_\Gamma=\Pi_{\Gamma^{(v,w)}}$.
\end{enumerate}
\end{lemma}

\begin{proof}
Straightforward.
\end{proof}

If $v,w$ are two vertices with $v+w$ in the isotropic radical of $Q$, then we write $v\equiv w$.
The relation $\equiv$ is an equivalence relation.

By $d_\Gamma(v)$ we denote the degree of a vertex $v$ in the graph $\Gamma$.

\begin{lemma} \label{valency}
Let $\Gamma$ be a finite connected  graph.
Let $v$ be a vertex of $\Gamma$ and $C: v= v_1,v_2,\ldots, v_n$ a minimal cycle on the distinct vertices $v_1,\dots,v_n$ in $\Gamma$.
Let $w_1,\ldots,w_l$ be the vertices outside $C$.

Then $$\Delta = \Gamma^{(v_2,v)\cdots(v_{n-1},\tau_{v_{n-2}}\cdots \tau_{v_2}(v))
(v_2,w_1)\cdots(v_{n-1},\tau_{v_{n-2}}\cdots \tau_{v_2}(w_1))\cdots
(v_2,w_l)\cdots(v_{n-1},\tau_{v_{n-2}}\cdots \tau_{v_2}(w_l))}$$ satisfies
$$d_\Delta(\tau_{v_{n-2}}\cdots \tau_{v_2}(v)) =d_\Gamma(v) -1.$$
\end{lemma}

\begin{proof}
Consider $$\Delta = \Gamma^{(v_2,v)\cdots(v_{n-1},\tau_{v_{n-2}}\cdots \tau_{v_2}(v))
(v_2,w_1)\cdots(v_{n-1},\tau_{v_{n-2}}\cdots \tau_{v_2}(w_1))\cdots
(v_2,w_l)\cdots(v_{n-1},\tau_{v_{n-2}}\cdots \tau_{v_2}(w_l))}$$

Under this transformation the vertex $v$ is changed to $v'=v_1+v_2+\ldots+v_{n-1}$.
All the other vertices of the cycle remain the same.
Now $f(v',v_k)=0$ for $2<k<n-1$ and $f(v',v_n)=0$. In particular, we find a tree of type $D_n$.

A vertex $w\not\in C$ is replaced by $$w'=w+f(w,v_1)v_2+f(w,v_1+v_2)v_2+\ldots+f(w,v_1+\ldots v_{n-1})v_{n-1}.$$
So, $$\begin{array}{ll}f(v_i,w')&=f(v_i,w+f(w,v_1)v_1+f(w,v_1+v_2)v_2+\ldots+f(w,v_1+\ldots v_{n-1})v_{n-1})\\
&=f(v_i,w)+f(w,v_1+\ldots+v_{i-1})+f(w,v_1+\ldots+v_{i+1})\\
&=f(v_{i+1},w),\\
\end{array}$$
for $2\leq i<n-1$, while 
$$f(v_{n-1},w')=f(v_1+\ldots+v_{n-1},w).$$
We find
$$\begin{array}{ll}
f(v',w')&=f(v_1+v_2+\ldots+v_{n-1},w')\\
&=f(v_1,w')+\ldots +f(v_{n-1},w')\\
&=f(v_2,w)+f(v_3,w)+\dots +f(v_{n-1},w)+f(v_1+\ldots+v_{n-1})\\
&=f(v_1,w)\\
&=f(v,w).\\
\end{array}$$
So, the degree of $v'$ in $\Delta$ is the same as the one for $v$ in $\Gamma$
except that the degree
 inside the transforms of the cycle $C$ has dropped from $2$ to $1$.
\end{proof}

\begin{theorem}\label{tree}
Every finite connected graph $\Gamma$ is equivalent to a tree
containing an  $E_6$ induced subgraph, or the quotient $\overline\Gamma$ is equivalent to $A_{n}$.  
\end{theorem}

\begin{proof}
Fix a maximal subset $T$ of $\mathcal{V}$ inducing a tree not containing two 
vertices of any cycle. Notice that, as a vertex is a tree not containing two 
vertices of a cycle,  such $T$ exist.
Let $v$ be
an endpoint of $T$ having a neighbor outside $T$. (If such $v$ does not exist, $\Gamma$ itself is a tree.)

Suppose $v$ lies on a cycle. Then, by the above, the cycle lies entirely in $(\mathcal{V}\setminus
T)\cup\{v\}$. Transformations involving vertices from $(\mathcal{V}\setminus T)\cup\{v\}$ have no  effect on the graph structure on $T$.
 So, we can apply  \cref{valency} to the induced graph on
$(\mathcal{V}\setminus T)\cup\{v\}$  to find a graph $\Delta$ equivalent  with 
$\Gamma$ by a transformation $\tau$  with $d_\Delta(\tau(v)) <d_\Gamma(v) $.

We replace $\Gamma$ by $\Delta$ and $v$ by $\tau(v)$. 
We can repeat this process only a finite number of times and 
end up in a situation where $v$ is not in a cycle.  We  pick a vertex
$w$ of $v$ outside $T$ and find the
induced graph on $T\cup\{w\}$ again to be a tree not containing two vertices of a cycle.  

Repeat the above
procedure, beginning with the choice of a vertex of $T$
separating the rest of $T$ from a nonempty remainder in $\Gamma$, until
there are no such end nodes. As at each step, the induced subtree becomes
bigger, this will terminate, and we have found a tree $T$ equivalent to $\Gamma$.

If the tree does not contain a subgraph $E_6$, then $T$ is of the form

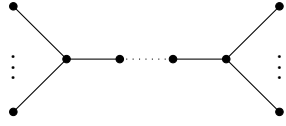
\begin{center}
\begin{tikzpicture}[scale=0.7]
\filldraw

(0,1) circle[radius=2pt]
(0,-1) circle[radius=2pt]
(1,0) circle[radius=2pt]
(2,0) circle[radius=2pt]
(3,0) circle[radius=2pt]
(4,0) circle[radius=2pt]
(5,1) circle[radius=2pt]
(5,-1) circle[radius=2pt];

\draw (0,1)--(1,0)--(0,-1);
\draw (1,0)--(2,0);
\draw (5,-1)--(4,0)--(5,1);
\draw[dotted] (2,0)--(3,0);
\draw (3,0)--(4,0);
\draw (0,0) node {$\vdots$};
\draw (5,0) node {$\vdots$};

\end{tikzpicture}
\end{center}
and, in $\overline T$, all end nodes (at the left or right) are mapped to a single vertex resulting into
a graph of type $A_n$. 
\end{proof}

\begin{theorem}\label{angraphthm}
A  finite connected graph $\Gamma$  is a line graph of a multi-graph if and only if $\overline\Gamma$ is equivalent to
$A_n$ for some integer $n$.
\end{theorem}

\begin{proof}
Suppose $\overline\Gamma$ 
is equivalent to $A_n$. Then $\overline{\Pi}_\Gamma=\Pi_{\overline\Gamma}$ is isomorphic to $\Pi_{A_n}$ and thus, by \cref{plsthm}, we find  $\Gamma$ to be a line graph.

If $\Gamma$ is equivalent to a tree containing and $E_6$ subdiagram, then
$\overline{\Pi}_\Gamma$ contains a subspace
isomorphic to $\Pi_{E_6}$, and hence it is not the cotriangular space of a set. See \cref{E6example}.
But then, again using  \cref{plsthm}, we find  $\overline\Gamma$ and also $\Gamma$ not to be line graphs.

By \cref{tree}, this proves the theorem.
\end{proof}

\section{Line graphs and $E_6$}
\label{sect:E6}

In this section we provide a proof for \cref{mainthm}.
We continue with the notation of previous sections.

\begin{lemma}\label{E6}
Let $\Gamma$ be a connected graph on $7$ points with $(V_\Gamma,Q_\Gamma)$ nondegenerate.
Then $\Gamma$ contains a subgraph equivalent to $E_6$,
or $\Gamma$ is equivalent to $A_7$.
\end{lemma}

\begin{proof}
If $\Gamma$ contains an induced subgraph $\Delta$ on $4$ points which is equivalent to $D_4$, then
the subspace $V_\Delta$ of $V_\Gamma$ contains $12$ anisotropic points and has  a 2-dimensional isotropic radical $R_\Delta$.
As the isotropic radical of $V_\Gamma$ is trivial, there are two vertices $v,w\in \mathcal{V}$ such that $\langle v,w\rangle^\perp\cap R=\{0\}$. (Here $\perp$ denotes orthogonality.)
The induced subgraph  $\Gamma_0$
on $\Delta$ and the vertices $v,w$ generates a nondegenerate hyperplane $V_{\Gamma_0}$ in which we find the subspace  $V_\Delta$.
This implies that 
 $(V_{\Gamma_0}, Q_{\Gamma_0})$ is $6$-dimensional and of $-$-type. See \cref{E6example}.

So, we can now assume that $\Gamma$ contains no subgraph equivalent to $D_4$.
Let $C$ be a maximal clique in $\Gamma$. 
A clique $C$ is equivalent to $A_{|C|}$.
So, we can assume $|C|\leq 6$. A vertex $v$ outside $C$ is adjacent to at most one  element
in $C$. For, otherwise, we find a subgraph
\begin{tikzpicture}[scale=0.2]
\filldraw[color=black] 
(0,0) circle [radius=2pt]
(1,0.5) circle [radius=2pt]
(1,-0.5) circle [radius=2pt]
(2,0) circle [radius=2pt];

\draw (0,0)--(1,0.5)--(2,0)--(1,-0.5)--(1,0.5);
\draw (0,0)--(1,-0.5);
\end{tikzpicture} equivalent to $D_4$.
Moreover, a point can not be in three maximal cliques, as this would mean that $\Gamma$ contains a subgraph
$D_4$.

Suppose $C_1,\dots, C_k$ are the maximal cliques of size at least $2$. 
As each vertex is in at most two such cliques, we find $k\leq 7$.
Assign to a vertex $v$ the pair $\{ i, j\}$, if $v$ is in the intersection $C_i\cap C_j$, where $i\neq j$.
If $v$ is in a unique maximal clique $C_i$, assign to it the pair  $\{i,l\}$ for some $k<l\leq 8$ such that different pairs are assigned to different vertices. This is possible by the above.
Two vertices are adjacent if and only if the assigned pairs meet nontrivially.
But then it is easy to check that $\Gamma$  is  equivalent to $A_7$.
\end{proof}

\begin{lemma}[Seven's Lemma, \cite{seven}]\label{seven}
If $\Gamma$ is a finite connected graph equivalent to a tree containing $E_6$, then $\Gamma$ contains a subgraph on $6$ vertices equivalent to $E_6$.
\end{lemma}

\begin{proof}
Suppose $\Gamma$ contains a subgraph $\Delta$ on $6$ vertices equivalent to $E_6$.
To prove the lemma it suffices to show that for any two vertices $v,w$ we will find $\Gamma^{(v,w)}$ also to contain a subgraph on $6$ vertices equivalent to $E_6$.

If $v,w$ are nonadjacent, or vertices both in $\Delta$ or both outside $\Delta$, this is clear.
So, consider the case where $w$ is a vertex of $\Delta$, adjacent to  $v$, which is not in  $\Delta$.
Then consider the graph $\widehat\Delta$ on the vertices of $\Delta^{(v,w)}$ and $w$.

If $(V_{\widehat\Delta},Q_{\widehat\Delta})$ a nondegenerate $7$-dimensional space, then, by \cref{E6}, we find that
the subgraph $\widehat\Delta$ of $\Gamma^{(v,w)}$ contains  a subgraph equivalent to $E_6$.
If $(V_{\widehat\Delta},Q_{\widehat\Delta})$ is a degenerate $7$-dimensional space, it contains a $1$-dimensional radical, and modulo this radical
we find a $6$-dimensional orthogonal space of $-$ type and we can take any $6$ vertices of $\widehat\Delta$ spanning a complement to the radical.
\end{proof}

\begin{theorem}\label{E6thm}
Let $\Gamma$ be a  connected graph.
Then 
$\Gamma$ is a line graph of a multi-graph, if and only if it does not contain an induced subgraph in $\mathcal{E}_6$.
\end{theorem}

\begin{proof}
If $\Gamma$ contains a subgraph in $\mathcal{E}_6$, it can not be a line graph.

Now assume that $\Gamma$ contains no  induced subgraph from $\mathcal{E}_6$.
Then $\overline\Gamma$ also has no induced subgraph from $\mathcal{E}_6$.
If $\overline\Gamma$ is a finite graph, then we can apply
\cref{tree} and \cref{seven}, and find that $\overline\Gamma$ is
equivalent to $A_n$ for some $n$ and hence, by \cref{angraphthm}, that $\Gamma$ is a line graph.

The case that $\Gamma$ is infinite can be dealt with in the following way.

Consider $\overline\Gamma$ which also has no induced subgraph from $\mathcal{E}_6$.
By \cref{plsthm} we only have to show that $\Pi_{\overline\Gamma}$ is (isomorphic to) a
cotriangular space of a set $\Omega$.
But by \cref{plsthm} this can be checked inside finite subspaces of $\Pi_{\overline\Gamma}$ and thus
inside finite subgraphs of $\overline\Gamma$.
\end{proof}

\cref{mainthm} follows from  \cref{E6thm}, as the graphs in \cref{E6graphs} are the graphs in the set $\mathcal{E}_6$.

\begin{remark}
To compute the list of graphs in $\mathcal{E}_6$ we have used the following observations.

If $\Gamma$ is an element of $\mathcal{E}_6$, then it does not contain any two nonadjacent vertices having the same neighbors, for otherwise, the sum of these two vertices is a nonzero vector in the radical of $Q_{\Gamma}$. Moreover, if $C$ is a minimal cycle in $\Gamma$, then there has to be a vertex that has an odd number of neighbors in the cycle, for otherwise the sum of the vertices in the cycle is in the radical of $Q_{\Gamma}$.
Finally, as follows from the second part of the  proof of \cref{E6}, the graph $\Gamma$ has an induced subgraph equivalent to $D_4$, so 
\begin{tikzpicture}[scale=0.2]
\filldraw[color=black] 
(0,0) circle [radius=2pt]
(1,0.5) circle [radius=2pt]
(1,-0.5) circle [radius=2pt]
(2,0) circle [radius=2pt];

\draw (0,0)--(1,0.5)--(2,0)--(1,-0.5)--(1,0.5);
\draw (0,0)--(1,-0.5);
\end{tikzpicture} or 
\begin{tikzpicture}[scale=0.2]
\filldraw[color=black] 
(0,0) circle [radius=2pt]
(1,1) circle [radius=2pt]
(1,0) circle [radius=2pt]
(2,0) circle [radius=2pt];

\draw (0,0)--(2,0);
\draw (1,0)--(1,1);
\end{tikzpicture}.

Starting with one of these two graphs on four vertices, one can add two extra vertices
to obtain a graph $\Gamma$ that satisfies the above requirements.
We have checked our list against the list of all connected graphs on $6$ vertices as given in \cite{6vertices}.
\end{remark}

\section{Line graphs of ordinary graphs}
\label{sect:linegraph}

In this section we determine the minimal graphs which are not line graphs of ordinary graphs and reprove Beineke's result \cref{beineke}. So, let $\Gamma$ be such a graph.
Assume that $\Gamma$ is $\mathcal{E}_6$-free, and hence by \cref{E6thm}   a  
line graph of a multi-graph $\Delta$.

As $\Gamma$ is not the line graph of an ordinary graph, we can assume that $\Delta$ contains two vertices, $v,w$ say, on which there are two or more edges.

If $\Delta$ contains a 3-fold edge, it contains
(by connectedness) a subgraph (not necessarily induced)
\begin{tikzpicture}[scale=0.5]
\filldraw
(0,0) circle[radius=2pt]
(1,0) circle[radius=2pt]
(2,0) circle[radius=2pt];

\draw (0,0)--(2,0);

\draw (1,0.05)--(2,0.05);
\draw (1,-0.05)--(2,-0.05);

\end{tikzpicture}
and we find $\Gamma$ to contain an induced subgraph $H_1$. So, we can assume that on each pair of points
there are at most two edges.

If $\Delta$ contains more than $4$ vertices,
then  we find at least one of the following three  subgraph:

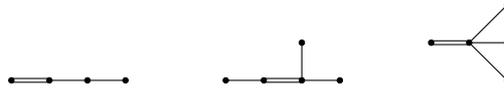
\begin{center}
\begin{tikzpicture}[scale=0.5]
\filldraw 
(0,0) circle[radius=2pt]
(1,0) circle[radius=2pt]
(2,0) circle[radius=2pt]
(3,0) circle[radius=2pt];

\draw (1,0)--(3,0);
\draw (0,0.05)--(1,0.05);
\draw (0,-0.05)--(1,-0.05);
\end{tikzpicture}
\hspace{1cm}
\begin{tikzpicture}[scale=0.5]
\filldraw 
(0,0) circle[radius=2pt]
(1,0) circle[radius=2pt]
(2,0) circle[radius=2pt]
(3,0) circle[radius=2pt]
(2,1) circle[radius=2pt];

\draw (0,0)--(1,0);
\draw (2,0)--(3,0);
\draw (2,0)--(2,1);
\draw (1,0.05)--(2,0.05);
\draw (1,-0.05)--(2,-0.05);
\end{tikzpicture}\hspace{1cm}
\begin{tikzpicture}[scale=0.5]
\filldraw 
(0,0) circle[radius=2pt]
(1,0) circle[radius=2pt]
(2,1) circle[radius=2pt]
(2,0) circle[radius=2pt]
(2,-1) circle[radius=2pt];

\draw (1,0)--(2,0);
\draw (1,0)--(2,1);
\draw (1,0)--(2,-1);

\draw (0,0.05)--(1,0.05);
\draw (0,-0.05)--(1,-0.05);
\end{tikzpicture}
\end{center}

But that implies that $\Gamma$
contains an induced subgraph $H_1$, $H_2$ or $H_3$.

If $\Delta$ contains $3$ or $4$ points it is straightforward to check that $\Gamma$ is not a line graph if and only if it does not contain induced subgraph $H_1$, $H_2$ or $H_3$.
This leads to the following theorem.

\begin{figure}

\begin{tikzpicture}[scale=0.5]
\filldraw[black] 

(0,0) circle [radius=2pt]
(1,0) circle [radius=2pt]
(2,0) circle [radius=2pt]
(1,1) circle [radius=2pt];

\draw (0,0)--(2,0);
\draw (1,0)--(1,1);

\draw (1,-2) node {$H_1=D_4$};

\filldraw[black] 

(4,0) circle [radius=2pt]
(5,0) circle [radius=2pt]
(6,0) circle [radius=2pt]
(5,1) circle [radius=2pt]
(5,-1) circle [radius=2pt];

\draw (4,0)--(5,1)--(6,0)--(5,-1)--(4,0);
\draw (4,0)--(6,0);
\draw (5,0)--(5,1);

\draw (5,-2) node {$H_2$};

\filldraw[black] 

(7,-1) circle [radius=2pt]
(9,-1) circle [radius=2pt]
(8,0) circle [radius=2pt]
(8,1) circle [radius=2pt]
(8,2) circle [radius=2pt];

\draw (7,-1)--(9,-1)--(8,2)--(7,-1)--(8,1)--(9,-1)--(8,0)--(7,-1);
\draw (8,0)--(8,2);

\draw (8,-2) node {$H_3$};

\end{tikzpicture}

\caption{Some forbidden subgraphs of a line graph.}
\label{forbiddenlinegraphs}
\end{figure}
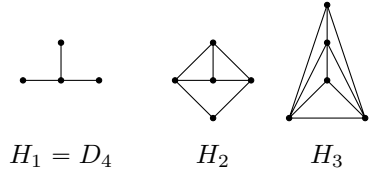

\begin{theorem}[Beineke, \cite{beineke}]\label{beineke}
A graph $\Gamma$ is an ordinary line graph if and only if it does not contain any of the graphs $H_1,H_2, H_3$ as in \cref{forbiddenlinegraphs},
or of the graphs $E_6^{(8)},E_6^{(12)},E_6^{(20)},E_6^{(22)}, E_6^{(25)}$ or $E_6^{(30)}$ from \cref{E6graphs}.
\end{theorem}

\begin{proof}
If a graph $\Gamma$ is the line graph of an ordinary graph, then clearly it does not contain any of the graphs $H_1,H_2, H_3$ as in \cref{forbiddenlinegraphs},
nor one of the graphs $E_6^{(8)},E_6^{(12)},E_6^{(20)},E_6^{(22)}, E_6^{(25)}$ or $E_6^{(30)}$ from \cref{E6graphs}.

The above shows that a graph $\Gamma$ which is not 
a line graph of an ordinary graph and does not contain  one of the graphs from $\mathcal{E}_6$, will contain $H_1,H_2$, or $H_3$.

As the graphs $E_6^{(8)},E_6^{(12)},E_6^{(20)},E_6^{(22)}, E_6^{(25)}$ or $E_6^{(30)}$ are precisely the graphs of $\mathcal{E}_6$ not containing an induced subgraph $H_1$, $H_2$ or $H_3$, we have proved the theorem.
\end{proof}

\section{Generalized line graphs}
\label{sect:genlinegraph}

In this section we consider generalized line graphs
and determine the minimal graphs which are not generalized line graphs.

\begin{definition}
A \emph{ cocktail party graph} is a graph in which every vertex is adjacent to precisely all, or  all but one other vertices.

Suppose $\Delta=(\mathcal{V},\mathcal{E})$ is a graph and for every vertex $v$ of $\Delta$ let $\Delta_v$ be a (possibly empty)  cocktail party graph.

Then the \emph{generalized line graph} obtained
from $\Delta$ and the various $\Delta_v$ for $v\in \mathcal{V}$
is the graph which is the union of the line graph
$L(\Delta)$ and all the graphs $\Delta_v$,
to which the following edges are  added:
a vertex $e$ of $L(\Delta)$ is adjacent to 
all vertices of $\Delta_v$, where $v$ is a vertex on $e$.
\end{definition}

This definition of generalized line graph
does not directly fit our purposes.

We use the following proposition, which shows that
 generalized line graphs are a special type of line graph of a multi-graph.

\begin{proposition}\label{generalizedgraphdef}
A graph $\Gamma$ is a generalized line graph if and only if it is the line graph of a 
multi-graph $\Delta$ such that
\begin{enumerate}
\item two vertices of $\Delta$ are on at most two common edges.
\item If $v,w$ are two vertices of $\Delta$ which are on two common edges, then
one of the vertices is on no other edges.
\end{enumerate}
\end{proposition}

\begin{proof}
Suppose we are given a generalized line graph $\Gamma$ 
based on a graph $\Delta=(\mathcal{V},\mathcal{E})$ and a collection
$\Delta_v$, where $v\in\mathcal{V}$, of  cocktail party graphs.

Then consider the multi-graph $\hat\Delta$
which consists of $\Delta$ to which 
we add  new vertices $\delta_{v,\{w,w'\}}$ for each $v\in \mathcal{V}$ and $\{w,w'\}$ a pair of nonadjacent vertices 
of $\Delta_v$, and 
two edges on $v$ and $\delta_{v,\{w,w'\}}$.
The graph $\Delta'$ satisfies (i) and (ii) and $\Gamma$ is the line graph of $\Delta'$.

On the other hand if $\Gamma$ is the line graph
of a graph $\hat\Delta$ satisfying (i) and (ii), then
let $\mathcal{V}_0$ be the vertices that are on a multiple
edge but not on another edge, and let $\Delta$ be the graph obtained
from $\hat\Delta$ by removing all vertices of $\mathcal{V}_0$ and the edges on them.
For each $v_0\in \mathcal{V}_0$, let $v_1$ be the 
unique neighbor. Then for each such $v_1$ define $\Delta_{v_1}$ to be the cocktail party graph
which is the line graph of the subgraph of $\hat\Delta$
induced on $v_1$ and all the $v_0\in \mathcal{V}_0$
adjacent to $v_1$. Then $\Gamma$ is the generalized line graph obtained from $\Delta$ and all these  cocktail party graphs.
\end{proof}

Now assume that $\Gamma$ is a connected graph which is not a generalized line graph, but $\mathcal{E}_6$-free.

We will prove that $\Gamma$ contains an induced subgraph isomorphic to $G_1,\dots, G_{10}$ or $G_{11}$.

By \cref{E6thm} we find that $\Gamma$ is the line graph of a connected multi-graph $\Delta$.
As we assume $\Gamma$ not to be a generalized line graph,
$\Delta$ will have to contain one of the following subgraphs:

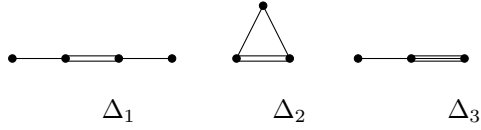
\begin{center}

\begin{tikzpicture}[scale=0.7] \filldraw 
(0,0) circle[radius=2pt]
(1,0) circle[radius=2pt]
(2,0) circle[radius=2pt]
(3,0) circle[radius=2pt];

\draw (0,0)--(1,0);
\draw (1,0.05)--(2,0.05);
\draw (1,-0.05)--(2,-0.05);

\draw (2,0)--(3,0);

\draw (2,-1) node {$\Delta_1$};
\end{tikzpicture} \hspace{0.5cm}
\begin{tikzpicture}[scale=0.7] \filldraw 
(1.5,1) circle[radius=2pt]
(1,0) circle[radius=2pt]
(2,0) circle[radius=2pt];

\draw (2,0)--(1.5,1)--(1,0);
\draw (1,0.05)--(2,0.05);
\draw (1,-0.05)--(2,-0.05);
\draw (2,-1) node {$\Delta_2$};
\end{tikzpicture}\hspace{0.5cm}\begin{tikzpicture}[scale=0.7] \filldraw 
(2,0) circle[radius=2pt]
(1,0) circle[radius=2pt]
(0,0) circle[radius=2pt];

\draw (1,0.05)--(2,0.05);
\draw (1,-0.05)--(2,-0.05);

\draw (0,0)--(2,0);
\draw (2,-1) node {$\Delta_3$};
\end{tikzpicture}
\end{center}

Suppose $\Delta$ contains a $5$-fold edge.
Then it contains a subgraph
\begin{tikzpicture}[scale=0.7] \filldraw 
(2,0) circle[radius=2pt]
(1,0) circle[radius=2pt]
(0,0) circle[radius=2pt];

\draw (1,0.05)--(2,0.05);
\draw (1,-0.04)--(2,-0.04);
\draw (1,0.08)--(2,0.08);
\draw (1,-0.08)--(2,-0.08);
\draw (0,0)--(2,0);
\end{tikzpicture}
and $\Gamma$ contains a subgraph $G_5$.

So, we can assume that $\Delta$ contains
only $k$-fold edges, where $k\leq 4$.

Now assume that $\Delta$ contains at least $5$ vertices.
If $\Delta$ contains $\Delta_1$, we can extend 
$\Delta_1$ to one of the following graphs: 

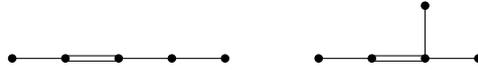
\begin{center}
\begin{tikzpicture}[scale=0.7]
\filldraw 
(0,0) circle[radius=2pt]
(1,0) circle[radius=2pt]
(2,0) circle[radius=2pt]
(3,0) circle[radius=2pt]
(4,0) circle[radius=2pt];

\draw (0,0)--(1,0);
\draw (2,0)--(4,0);
\draw (1,0.05)--(2,0.05);
\draw (1,-0.05)--(2,-0.05);
\end{tikzpicture}\hspace{1cm}
\begin{tikzpicture}[scale=0.7]
\filldraw 
(0,0) circle[radius=2pt]
(1,0) circle[radius=2pt]
(2,0) circle[radius=2pt]
(3,0) circle[radius=2pt]
(2,1) circle[radius=2pt];

\draw (0,0)--(1,0);
\draw (2,0)--(2,1);
\draw (2,0)--(3,0);
\draw (1,0.05)--(2,0.05);
\draw (1,-0.05)--(2,-0.05);
\end{tikzpicture}

\end{center}

But then $\Gamma$ contains $G_1$ or $G_3$.

So, from now on we can and do assume that $\Delta$ contains no $\Delta_1$.
If $\Delta$ contains $\Delta_3$, then by adding some edges,
we can extend it to one of the following subgraphs:
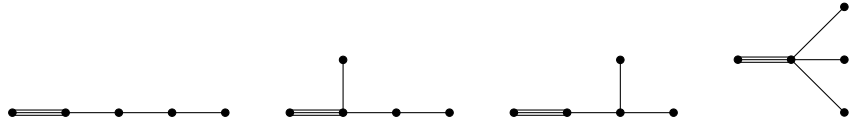
\begin{center}
\begin{tikzpicture}[scale=0.7]
\filldraw 
(0,0) circle[radius=2pt]
(1,0) circle[radius=2pt]
(2,0) circle[radius=2pt]
(3,0) circle[radius=2pt]
(4,0) circle[radius=2pt];

\draw (0,0)--(4,0);
\draw (0,0.05)--(1,0.05);
\draw (0,-0.05)--(1,-0.05);
\end{tikzpicture}
\hspace{0.5cm}
\begin{tikzpicture}[scale=0.7]
\filldraw 
(0,0) circle[radius=2pt]
(1,0) circle[radius=2pt]
(2,0) circle[radius=2pt]
(1,1) circle[radius=2pt]
(3,0) circle[radius=2pt];

\draw (0,0)--(3,0);
\draw (1,0)--(1,1);
\draw (0,0.05)--(1,0.05);
\draw (0,-0.05)--(1,-0.05);
\end{tikzpicture}
\hspace{0.5cm}
\begin{tikzpicture}[scale=0.7]
\filldraw 
(0,0) circle[radius=2pt]
(1,0) circle[radius=2pt]
(2,0) circle[radius=2pt]
(2,1) circle[radius=2pt]
(3,0) circle[radius=2pt];

\draw (0,0)--(3,0);
\draw (2,0)--(2,1);
\draw (0,0.05)--(1,0.05);
\draw (0,-0.05)--(1,-0.05);
\end{tikzpicture}
\hspace{0.5cm}
\begin{tikzpicture}[scale=0.7]
\filldraw 
(0,0) circle[radius=2pt]
(1,0) circle[radius=2pt]
(2,0) circle[radius=2pt]
(2,1) circle[radius=2pt]
(2,-1) circle[radius=2pt];

\draw (0,0)--(2,0);
\draw (1,0)--(2,1);
\draw (1,0)--(2,-1);
\draw (0,0.05)--(1,0.05);
\draw (0,-0.05)--(1,-0.05);
\end{tikzpicture}
\end{center}
But then we find in $\Gamma$ an induced subgraph  $G_4$, $G_8$, $G_6$, or $G_{11}$, respectively.

So, we can also assume that $\Delta$ does not contain
$\Delta_3$.
It remains to consider the case that $\Delta$ contains $\Delta_2$.

In this case we can extend $\Delta_2$ to one of the following subgraphs
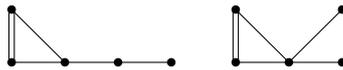
\begin{center}
\begin{tikzpicture}[scale=0.7]
\filldraw 
(0,0) circle[radius=2pt]
(0,1) circle[radius=2pt]
(1,0) circle[radius=2pt]
(2,0) circle[radius=2pt]
(3,0) circle[radius=2pt];

\draw (0.05,0)--(0.05,1);
\draw (-0.05,0)--(-0.05,1);
\draw (0,1)--(1,0);
\draw (0,0)--(3,0);
\end{tikzpicture}
\hspace{0.5cm}
\begin{tikzpicture}[scale=0.7]
\filldraw 
(0,0) circle[radius=2pt]
(0,1) circle[radius=2pt]
(1,0) circle[radius=2pt]
(2,0) circle[radius=2pt]
(2,1) circle[radius=2pt];

\draw (0.05,0)--(0.05,1);
\draw (-0.05,0)--(-0.05,1);
\draw (0,1)--(1,0);
\draw (0,0)--(2,0);
\draw (1,0)--(2,1);
\end{tikzpicture}
\hspace{0.5cm}

\end{center}
and find subgraphs $G_7$ or $G_{10}$ inside $\Gamma$.

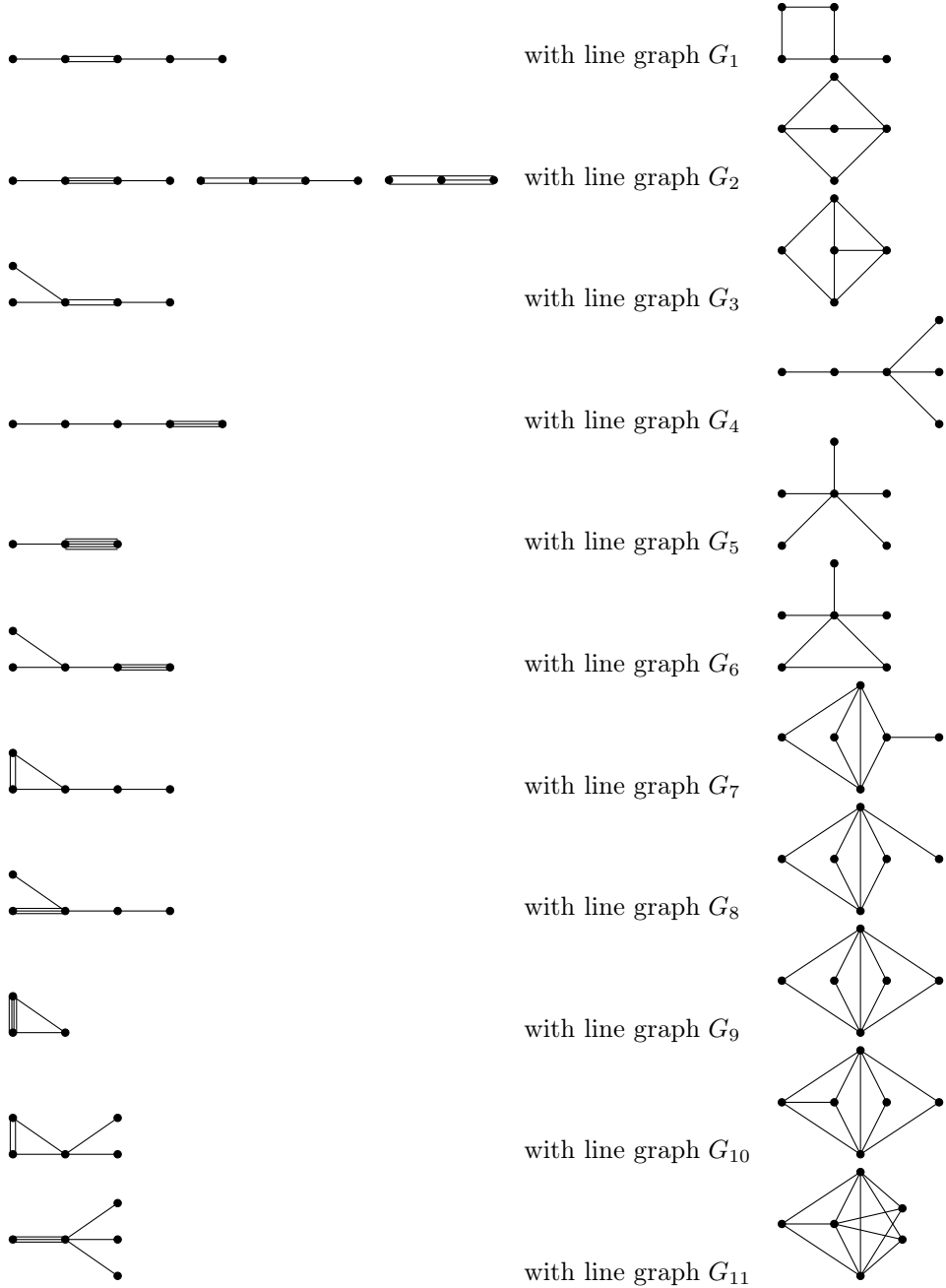
\begin{figure}
\begin{tabular}{lll}
\begin{tikzpicture}[scale=0.7]
\filldraw 
(0,0) circle[radius=2pt]
(1,0) circle[radius=2pt]
(2,0) circle[radius=2pt]
(3,0) circle[radius=2pt]
(4,0) circle[radius=2pt];

\draw (0,0)--(1,0);
\draw (2,0)--(4,0);
\draw (1,0.05)--(2,0.05);
\draw (1,-0.05)--(2,-0.05);
\end{tikzpicture}& with line graph $G_1$&\begin{tikzpicture}[scale=0.7]
\filldraw 
(0,0) circle[radius=2pt]
(1,0) circle[radius=2pt]
(2,0) circle[radius=2pt]
(1,1) circle[radius=2pt]
(0,1) circle[radius=2pt];

\draw (0,0)--(1,0)--(1,1)--(0,1)--(0,0);
\draw (1,0)--(2,0);
\end{tikzpicture}\\

\begin{tikzpicture}[scale=0.7]
\filldraw 
(0,0) circle[radius=2pt]
(1,0) circle[radius=2pt]
(2,0) circle[radius=2pt]
(3,0) circle[radius=2pt];

\draw (0,0)--(3,0);
\draw (1,0.05)--(2,0.05);
\draw (1,-0.05)--(2,-0.05);
\end{tikzpicture}\hspace{0.3cm}\begin{tikzpicture}[scale=0.7]
\filldraw 
(0,0) circle[radius=2pt]
(1,0) circle[radius=2pt]
(2,0) circle[radius=2pt]
(3,0) circle[radius=2pt];

\draw (2,0)--(3,0);
\draw (0,0.05)--(2,0.05);
\draw (0,-0.05)--(2,-0.05);
\end{tikzpicture}\hspace{0.3cm}\begin{tikzpicture}[scale=0.7]
\filldraw 
(0,0) circle[radius=2pt]
(1,0) circle[radius=2pt]
(2,0) circle[radius=2pt];

\draw (0,0.08)--(2,0.08);
\draw (0,-0.08)--(2,-0.08);
\draw (1,0)--(2,0);
\end{tikzpicture}
&
with line graph $G_2$&\begin{tikzpicture}[scale=0.7]
\filldraw 
(0,1) circle[radius=2pt]
(1,0) circle[radius=2pt]
(2,1) circle[radius=2pt]
(1,1) circle[radius=2pt]
(1,2) circle[radius=2pt];

\draw (0,1)--(1,0)--(2,1)--(1,2)--(0,1);
\draw (0,1)--(1,1)--(2,1);
\end{tikzpicture}
\\
\begin{tikzpicture}[scale=0.7]
\filldraw 
(0,0) circle[radius=2pt]
(1,0) circle[radius=2pt]
(2,0) circle[radius=2pt]
(3,0) circle[radius=2pt]
(0,0.7) circle[radius=2pt];

\draw (0,0)--(1,0);
\draw (2,0)--(3,0);
\draw (1,0.05)--(2,0.05);
\draw (1,-0.05)--(2,-0.05);
\draw (0,0.7)--(1,0);

\end{tikzpicture}&
with line graph $G_3$&\begin{tikzpicture}[scale=0.7]
\filldraw 
(0,1) circle[radius=2pt]
(1,0) circle[radius=2pt]
(2,1) circle[radius=2pt]
(1,1) circle[radius=2pt]
(1,2) circle[radius=2pt];

\draw (0,1)--(1,0)--(2,1)--(1,2)--(0,1);
\draw (1,2)--(1,0);
\draw  (1,1)--(2,1);
\end{tikzpicture}\\
\begin{tikzpicture}[scale=0.7]
\filldraw 
(0,0) circle[radius=2pt]
(1,0) circle[radius=2pt]
(2,0) circle[radius=2pt]
(3,0) circle[radius=2pt]
(4,0) circle[radius=2pt];

\draw (0,0)--(4,0);
\draw (3,0.05)--(4,0.05);
\draw (3,-0.05)--(4,-0.05);
\end{tikzpicture}&
with line graph $G_4$&\begin{tikzpicture}[scale=0.7]
\filldraw 
(1,1) circle[radius=2pt]
(2,1) circle[radius=2pt]
(3,1) circle[radius=2pt]
(4,1) circle[radius=2pt]
(4,2) circle[radius=2pt]
(4,0) circle[radius=2pt];

\draw (1,1)--(4,1);
\draw (3,1)--(4,2);
\draw (3,1)--(4,0);

\end{tikzpicture}\\
\begin{tikzpicture}[scale=0.7]
\filldraw 
(0,0) circle[radius=2pt]
(1,0) circle[radius=2pt]
(2,0) circle[radius=2pt];

\draw (0,0)--(2,0);
\draw (1,0.05)--(2,0.05);
\draw (1,-0.05)--(2,-0.05);
\draw (1,0.1)--(2,0.1);
\draw (1,-0.1)--(2,-.1);
\end{tikzpicture}
&
with line graph $G_5$&\begin{tikzpicture}[scale=0.7]
\filldraw 
(0,0) circle[radius=2pt]
(1,1) circle[radius=2pt]
(2,0) circle[radius=2pt]
(0,1) circle[radius=2pt]
(2,1) circle[radius=2pt]
(1,2) circle[radius=2pt];

\draw (0,0)--(1,1);
\draw (2,0)--(1,1);
\draw (2,1)--(1,1);
\draw (0,1)--(1,1);
\draw (1,2)--(1,1);
\end{tikzpicture}\\
\begin{tikzpicture}[scale=0.7]
\filldraw 
(0,0) circle[radius=2pt]
(1,0) circle[radius=2pt]
(2,0) circle[radius=2pt]
(3,0) circle[radius=2pt]
(0,0.7) circle[radius=2pt];

\draw (0,0)--(3,0);
\draw (2,0.05)--(3,0.05);
\draw (2,-0.05)--(3,-0.05);
\draw (1,0)--(0,0.7);
\end{tikzpicture}
&
with line graph $G_6$&\begin{tikzpicture}[scale=0.7]
\filldraw 
(0,0) circle[radius=2pt]
(1,1) circle[radius=2pt]
(2,0) circle[radius=2pt]
(0,1) circle[radius=2pt]
(2,1) circle[radius=2pt]
(1,2) circle[radius=2pt];

\draw (0,0)--(1,1);
\draw (2,0)--(1,1);
\draw (2,1)--(1,1);
\draw (0,1)--(1,1);
\draw (1,2)--(1,1);
\draw (0,0)--(2,0);
\end{tikzpicture}\\
\begin{tikzpicture}[scale=0.7]
\filldraw 
(0,0) circle[radius=2pt]
(1,0) circle[radius=2pt]
(2,0) circle[radius=2pt]
(3,0) circle[radius=2pt]
(0,0.7) circle[radius=2pt];

\draw (0,0)--(3,0);
\draw(0,0.7)--(1,0);
\draw (0.05,0)--(0.05,0.7);
\draw (-0.05,0)--(-0.05,0.7);
\end{tikzpicture}
&
with line graph $G_7$&\begin{tikzpicture}[scale=0.7]
\filldraw 
(0,1) circle[radius=2pt]
(1,1) circle[radius=2pt]
(1.5,0) circle[radius=2pt]
(1.5,2) circle[radius=2pt]
(2,1) circle[radius=2pt]
(3,1) circle[radius=2pt];

\draw (0,1)--(1.5,0)--(1.5,2)--(0,1);
\draw (1.5,0)--(1,1)--(1.5,2);
\draw (1.5,0)--(2,1)--(1.5,2);
\draw (2,1)--(3,1);
\end{tikzpicture}\\
\begin{tikzpicture}[scale=0.7]
\filldraw 
(0,0) circle[radius=2pt]
(1,0) circle[radius=2pt]
(2,0) circle[radius=2pt]
(3,0) circle[radius=2pt]
(0,0.7) circle[radius=2pt];

\draw (0,0)--(3,0);
\draw(0,0.7)--(1,0);
\draw (0,0.05)--(1,0.05);
\draw (0,-0.05)--(1,-0.05);

\end{tikzpicture}
&
with line graph $G_8$&\begin{tikzpicture}[scale=0.7]
\filldraw 
(0,1) circle[radius=2pt]
(1,1) circle[radius=2pt]
(1.5,0) circle[radius=2pt]
(1.5,2) circle[radius=2pt]
(2,1) circle[radius=2pt]
(3,1) circle[radius=2pt];

\draw (0,1)--(1.5,0)--(1.5,2)--(0,1);
\draw (1.5,0)--(1,1)--(1.5,2);
\draw (1.5,0)--(2,1)--(1.5,2);
\draw (1.5,2)--(3,1);
\end{tikzpicture}\\
\begin{tikzpicture}[scale=0.7]
\filldraw 
(0,0) circle[radius=2pt]
(1,0) circle[radius=2pt]
(0,0.7) circle[radius=2pt];

\draw (0,0)--(1,0);
\draw(0,0.7)--(1,0);

\draw (-0.07,0)--(-0.07,0.7);
\draw (0.07,0)--(0.07,0.7);
\draw (-0.02,0)--(-0.02,0.7);
\draw (0.02,0)--(0.02,0.7);

\end{tikzpicture}
&
with line graph $G_9$&\begin{tikzpicture}[scale=0.7]
\filldraw 
(0,1) circle[radius=2pt]
(1,1) circle[radius=2pt]
(1.5,0) circle[radius=2pt]
(1.5,2) circle[radius=2pt]
(2,1) circle[radius=2pt]
(3,1) circle[radius=2pt];

\draw (0,1)--(1.5,0)--(1.5,2)--(0,1);
\draw (1.5,0)--(1,1)--(1.5,2);
\draw (1.5,0)--(2,1)--(1.5,2);
\draw (1.5,2)--(3,1)--(1.5,0);
\end{tikzpicture}\\
\begin{tikzpicture}[scale=0.7]
\filldraw 
(0,0) circle[radius=2pt]
(1,0) circle[radius=2pt]
(2,0) circle[radius=2pt]
(2,0.7) circle[radius=2pt]
(0,0.7) circle[radius=2pt];

\draw (0,0)--(2,0);
\draw (1,0)--(2,0.7);
\draw(0,0.7)--(1,0);
\draw (0.05,0)--(0.05,0.7);
\draw (-0.05,0)--(-0.05,0.7);

\end{tikzpicture}
&
with line graph $G_{10}$&\begin{tikzpicture}[scale=0.7]
\filldraw 
(0,1) circle[radius=2pt]
(1,1) circle[radius=2pt]
(1.5,0) circle[radius=2pt]
(1.5,2) circle[radius=2pt]
(2,1) circle[radius=2pt]
(3,1) circle[radius=2pt];

\draw (0,1)--(1.5,0)--(1.5,2)--(0,1);
\draw (1.5,0)--(1,1)--(1.5,2);
\draw (1.5,0)--(2,1)--(1.5,2);
\draw (1.5,2)--(3,1)--(1.5,0);
\draw (0,1)--(1,1);
\end{tikzpicture}\\
\begin{tikzpicture}[scale=0.7]
\filldraw 
(0,0) circle[radius=2pt]
(1,0) circle[radius=2pt]
(2,0) circle[radius=2pt]
(2,0.7) circle[radius=2pt]
(2,-0.7) circle[radius=2pt];

\draw (0,0)--(2,0);
\draw (1,0)--(2,0.7);
\draw (1,0)--(2,-0.7);
\draw (0,0.05)--(1,0.05);
\draw (0,-0.05)--(1,-0.05);

\end{tikzpicture}
&
with line graph $G_{11}$&\begin{tikzpicture}[scale=0.7]
\filldraw 
(0,1) circle[radius=2pt]
(1,1) circle[radius=2pt]
(1.5,0) circle[radius=2pt]
(1.5,2) circle[radius=2pt]
(2.3,1.3) circle[radius=2pt]
(2.3,0.7) circle[radius=2pt];

\draw (0,1)--(1.5,0)--(1.5,2)--(0,1);
\draw (1.5,0)--(1,1)--(1.5,2);
\draw (1.5,0)--(2.3,1.3)--(1.5,2);
\draw (1.5,2)--(2.3,0.7)--(1.5,0);
\draw (0,1)--(1,1);
\draw (1,1)--(2.3,0.7);
\draw (1,1)--(2.3,1.3);
\end{tikzpicture}\\

\end{tabular}
\caption{Some forbidden graphs $\Delta$ and their line graphs $G_1,\dots,G_{11}$.}
\label{deltagraphs}\label{Ggraphs}
\end{figure}

If $\Delta$ contains $\leq 4$ vertices and no
$k$-fold edges for $k>4$, then we can easily
check that $\Gamma$ is a generalized line graph if and only if it does not contain any induced subgraph
isomorphic to $G_1,\dots, G_{10}$ or $G_{11}$.
For some examples, see \cref{Ggraphs}.

We are now in a position to prove the following:

\begin{theorem}[\cite{GLG}]\label{GLGthm1}
A connected graph $\Gamma$ is a generalized line
graph if and only if it does not contain one of the $31$ graphs
$G_1,\dots, G_{11}$ or $E_6^{(i)}$ of \cref{E6graphs}, where $i$ is an element of $\{1,2,4,8,9,10,12,16,17,20,21,23,24,25,27,28,30,31,32,33\}$.
\end{theorem}

\begin{proof}
By the above we find that $\Gamma$ is a generalized line
graph if and only if it does not contain one of the graphs from $\mathcal{E}_6$ or $G_1,\dots,G_{11}$.
Removing from $\mathcal{E}_6$ graphs that contain one of $G_1,\dots,G_{11}$ leaves us with the additional $20$ graphs $E_6^{(i)}$ of \cref{E6graphs}, where $i\in \{1,2,4,8,9,10,12,16,17,20,21,23,24,25,27,28,30,31,32,33\}$.
\end{proof}

\section{Concluding remarks}
\label{sect:remarks}

We end this paper with some remarks.

\begin{enumerate}
\item By \cref{mainthm} checking whether a finite graph is the line graph of a multi-graph, only requires to consider all subgraphs on $6$ vertices. Hence this can be done in polynomial time.
An algorithm to find a multi-graph $\Delta$ with prescribed line graph $L(\Delta)$ is given in \cite{cuypers_whitney}. Notice that $\Delta$ need not be uniquely determined.

\item
In \cite{vija},  Vijayakumar considers root lattices generated by  sets of roots and proves that 
a set of roots of equals length, which can not be embedded in a root system of type $D$, contains six roots generating the root lattice $E_6$. See \cite[Theorem 16]{vija}.
This result can easily be shown to be a consequence of Seven's \cref{seven}.

\item
The infimum of the least eigenvalues of all finite induced subgraphs of an infinite graph is defined to be its least eigenvalue. Finite generalized line graphs are graphs with least eigenvalue $\geq -2$.
By our results, together with the fact that any finite connected graph  with  least eigenvalue $\geq -2$ is a generalized line graph or can be represented by roots in the root system $E_8$ and has at most $36$ vertices, see \cite{css},
we immediately find that  a connected  infinite graph with least eigenvalue  $\geq -2$ is a generalized line graph (see \cite{vija2}). Indeed,
every finite subgraph of an infinite graph $\Gamma$ can be embedded into a connected subgraph on more than $36$ vertices and therefor is a generalized line graph and hence contains none of the $31$
graphs of \cref{GLGthm}. But then \cref{GLGthm}
also shows that $\Gamma$ itself is a generalized line graph.
 
\item In \cite{signed} Vijayakumar considers edge signed graphs that can be embedded in a (possibly infinite dimensional) root lattice of type $D$. He characterizes these graphs by a list of 49 forbidded subgraphs.
Theorem \ref{mainthm}
and methods similar to those used   in \cref{sect:linegraph} and \cref{sect:genlinegraph}
can be used to give a new proof for this characterization of these graphs. 

\item Given a finite simply laced Coxeter diagram $\Gamma$, we find, by \cref{tree}, that the corresponding Coxeter group $(W(\Gamma),S)$ admits a quotient that maps the generators in $S$ to transpositions
in a symmetric group if and only if the graph $\Gamma$ is a line graph of a multi-graph.
If the corresponding Coxeter group $(W(\Gamma),S)$ does not have such quotient, then there are six generators $s_1,\dots, s_6$ such that $\langle s_1,\dots,s_6\rangle$ admits a quotient $W(E_6)$
which maps the six generators to reflections in the Weyl group $W(E_6)$.

\end{enumerate}

 \bibliographystyle{plain}

\bibliography{rootgraphs.bib}

\vspace{2cm}

\parindent=0pt
Hans Cuypers\\
Department of Computer Science and Mathematics\\
Eindhoven University of Technology\\
P.O. Box 513 5600 MB, Eindhoven\\
The Netherlands\\
email: f.g.m.t.cuypers@tue.nl

 \end{document}